\documentclass{article}

\usepackage[latin1]{inputenc}
\usepackage[english]{babel}
\usepackage{amsfonts}
\usepackage{amssymb}
\usepackage{amsmath}
\usepackage{xypic}
\usepackage[all]{xy}
\usepackage{mathrsfs}
\usepackage{latexsym}
\usepackage{amscd}
\usepackage{theorem}
\usepackage{QEDeng}

\usepackage[mathscr]{eucal}

\usepackage{setspace}                    % Line spacing

\usepackage{tocbibind}

\usepackage{hyperref}

\newtheorem{prop}{Proposition}[section]
\newtheorem{lem}[prop]{Lemma}
\newtheorem{cor}[prop]{Corollary}
\newtheorem{thm}[prop]{Theorem}
\newtheorem{dfn}[prop]{Definition}

\theorembodyfont{\rm}
\newtheorem{rmk}[prop]{Remark}
\newcommand{\nn}{\mathbb{N}}
\newcommand{\zz}{\mathbb{Z}}

\newcommand{\pp}{\mathbb{P}}
\newcommand{\PP}{\mathbb{P}}
\newcommand{\cou}{{C}}

\newcommand{\picard}{\mathrm{Pic}}

\newcommand{\gr}{\mathbb{G}\mathrm{r}}

\newcommand{\oo}{\mathcal{O}}

\newcommand{\lin}{\mathcal{L}}

\newcommand{\iso}{\tilde{\to}}

\newcommand{\pic}{\mathrm{Pic}}

\newcommand{\Hom}{{\mathrm{Hom}}}

\newcommand{\de}{\mathrm{d}}

\newcommand{\rk}{\mathrm{rk}}
\newcommand{\cod}{\mathrm{codim}}

\newcommand{\turndown}[1]{\rotatebox[origin=c]{270}{\ensuremath#1}}

\newcommand{\hookdownarrow}{\turndown{\hookrightarrow}} 

         \newif\ifpdf
        \ifx\pdfoutput\undefined
        \pdffalse % we are not running PDFLaTeX
        \else
        \pdfoutput=1 % we are running PDFLaTeX
        \pdftrue
        \fi

        \ifpdf
        \usepackage[pdftex]{graphicx}
        \else
        \usepackage{graphicx}
        \fi

\newcommand{\grafica}{
        \ifpdf
        \DeclareGraphicsExtensions{.pdf, .jpg}
        \else
        \DeclareGraphicsExtensions{.eps, .jpg}
        \fi}
\title{Stability of line bundles transforms on curves with respect to
 low codimensional  subspaces }
\author{Ernesto Mistretta}
\date{}

\begin{document}

\grafica

\maketitle

% \section{Butler's Theorem}

\begin{abstract}

We show the stability of certain syzygies of line bundles on curves,
which we call transforms, and are kernels of the evaluation map on subspaces of the space of global sections. For the transforms constructed, we prove the existence of 
reducible theta divisors, in the cases where the slope is integer.

\end{abstract}

\section{Introduction}

In the study of vector bundles on curves, it is a natural question
to investigate the stability of kernels of evaluation maps of global
sections.
This was used in particular by 
Paranjape and Ramanan (cf. \cite{pr}), and
Butler (cf. \cite{but}), to prove normal
generation of certain vector bundles, by Ein and Lazarsfeld 
(cf. \cite{e-l}) to show
the stability of the Picard bundle, and by Beauville
(\emph{e.g.} in \cite{bove}) to study theta divisors.

 \begin{dfn}
 Let $\cou$ be a smooth projective curve 
 over an algebraically closed field $\Bbbk$, 
 and $E$ a globally generated vector bundle over $\cou$.
 We call 
 $M_{V,E} := \mathrm{ker}( V \otimes \oo_{\cou}
 \twoheadrightarrow E)$
 the \emph{transform} of the vector bundle $E$ 
 with respect to the generating subspace $V \subset H^0(\cou, E)$,
 and 
 $M_E := M_{H^0(E),E} = 
 \mathrm{ker}( H^0 (\cou, E) \otimes \oo_{\cou}
 \twoheadrightarrow E)$ 
 the \emph{total transform} of $E$.
 \end{dfn}

Starting from a result of Butler, 
who proved the stability of total transforms 
under certain hypothesis,
we want to investigate the stability of 
transforms of line bundles by generic subspaces of certain codimensions.

 \begin{thm}[Butler]
 \label{but}
 Let $\cou$ be a smooth projective curve of genus 
 $g \geqslant 1$ 
 over an algebraically closed field $\Bbbk$, 
 and $E$ a semistable vector bundle over $\cou$ 
 with slope $\mu (E) \geqslant 2g$,
 then the vector bundle 
 $M_E := \mathrm{ker}( H^0 (\cou, E) \otimes \oo_{\cou}
 \twoheadrightarrow E)$
 is semistable.
 Furthermore, if $E$ is stable and $\mu (E) \geqslant 2g$,
 then $M_E$ is stable, unless $\mu (E) = 2g$, 
 and either $C$ is hyperelliptic or 
 $\omega_C \hookrightarrow E$.
 \end{thm}

It is natural to ask what happens taking subspaces in the place 
of the vector space of global sections.
Our results can be resumed to the following theorem:

\begin{thm}

Let $\lin$ be a line bundle  of degree $d$ on a curve 
$\cou$ of genus $g \geqslant 2$, such that  
$  d \geqslant 2g + 2c$, with $1 \leqslant c \leqslant g$.
Then $M_{V, \lin}$ is semistable for a 
generic
subspace $V \subset H^0 (\lin)$ of codimension $c$. 
It is stable unless $d= 2g + 2c$ and the curve is hyperelliptic.

\end{thm}

Similar results can be deduced by some constructions in Vincent
Mercat's work \cite{mer} on Brill-Noether's loci, but we think that in our case
it is useful to give a more direct proof which applies to all line
bundles of degree $d \geqslant 2g+2c$ and not only generic ones.

Eventually,
we observe the existence of theta divisors associated to the
(semi)stable transforms having integer slope $-2$. 
Those theta divisors are always non integral, and in most cases reducible, 
hence give further examples of
stable vector bundles admitting a reducible theta divisor
(cf. \cite{bove}).

\begin{rmk}

A geometrical interpretation of those kinds of results
goes as follows:
a generating subspace $V \subset H^0(C,\lin)$ gives rise to 
a base point free linear system $|V| \subset |{\lin}|$ on the curve $C$,
and determines a map
$\varphi_V \colon C \to \pp(V^*)$, 
which asociates to a point $x \in C$
the hyperplane of global sections in $V$ vanishing in $x$.
The Euler sequence on $\pp(V^*)$ is the dual of the tautological
sequence:
$$ 0 \to \Omega_{\pp(V^*)} (1) \to V \otimes \oo_{\pp(V^*)} 
\to \oo_{\pp(V^*)}(1) \to 0$$
which restricted to $C$ gives the evaluation sequence
$$0 \to M_{V,{\lin}} \to V \otimes \oo_C \to {\lin} \to 0~.$$

As stability of a vector bundle is not affected 
by dualizing and tensorizing by a line bundle,
we see that stability of 
$M_{V,{\lin}} = \Omega (1)_{|C}$ is equivalent to the stability of the restriction 
of the tangent bundle of the projective space $\pp(V^*)$ to the curve $C$.

\end{rmk}

So our theorem translates to 

\begin{thm}
\label{tangrestr}

Let $C \subset \pp^{d-g}$ be a genus $g\geqslant 2$ degree $d$ 
non-degenerate smooth curve, 
where $d > 2g + 2c$, and $c$ is a constant such that $1 \leqslant c \leqslant g$.
Then 
for the generic projection $\pp^{d-g} \dashrightarrow \pp^{d-g-c}$
the restriction $T_{\PP^{d-g-c} |C}$ is stable.
\end{thm}

\subsection{Acknowledgments}
I wish to thank my advisor Daniel Huybrechts
for his help during my thesis,
Arnaud Beauville for his interest and his comments,
and my colleagues of Paris and Rome for their friendly advices.

\section{Stability of transforms}

We essentially use the following two lemmas:

\begin{lem}[Butler]
\label{butlem}
 Let $\cou$ be a curve of genus $g \geqslant 2$, 
 $F$ a vector bundle on $\cou$ with no trivial summands, 
 and such that $h^1(F) \neq 0$.
 Suppose that $V \subset H^0(F)$ generates $F$.
 If $N=M_{V,F}$ is stable, 
 then $\mu(N) \leqslant -2$.
 Furthermore, $\mu(N)=-2$ implies that either 
 $\cou$ is hyperelliptic  $F$ is the hyperelliptic bundle
 and $N$  its dual, 
 or $F = \omega$ and $N = M_{\omega}$. 
\end{lem}

The proof of this lemma
is based on the result by Paranjape, 
Ramanan asserting the stability of $M_{\omega}$
(see \cite{but} and \cite{pr}).

\begin{lem}
\label{part1}

Let $\lin$ be a degree $d \geqslant 2g+2c$ 
line bundle on a
curve $\cou$ of genus $g \geqslant 2$, with $c \leqslant g$
and let $V \subset H^0(\lin)$ be a generating subspace of codimension $c$.
Supppose there exists 
a stable subbundle of maximal slope
$N \hookrightarrow M_{V, \lin}$
such that
$0 \neq N \neq M_{V,\lin}$ and $\mu(N) \geqslant \mu(M_{V,\lin})$.

Then there exists a line bundle $F$ of degree 
$f \leqslant d-1$, a generating subspace $W \subset  H^0(F)$,
and an injection 
$F \hookrightarrow \lin$
such that $N$ fits into the following commutative diagram
$$
\begin{array}{ccccccccc}
~~~ 0 & \to & N & \to & W \otimes \oo_{\cou} & \to & F
& \to & 0 \\
  &  & \hookdownarrow &  & \hookdownarrow &
 & \hookdownarrow & &  \\
 0 & \to & M_{V, \lin} & \to & V \otimes \oo_{\cou} &
\to & \lin & \to & 0 ~,
\end{array}
$$
\emph{i.e.} a destabilization of $M_{V, \lin}$ must be
the transform of a line bundle 
injecting into $\lin$ such that the global sections 
we are transforming by are in $V$.

\end{lem}

The importance of this lemma lies in the fact that we associate 
a line bundle $F$ to a destabilizing $N$, 
and this allows us more easily to parametrize destabilizations and bound their 
dimension.

\begin{proof}
We remark  that $\mu( M_{V, \lin}) = -d / (d-g-c) \geqslant -2$ 
for $d \geqslant 2g +2c$.
Consider a stable subbundle 
$N \hookrightarrow M_{V, \lin}$ of maximal slope.
Then it fits into the commutative diagram
$$
\begin{array}{ccccccccc}
~~~ 0 & \to & N & \to & W \otimes \oo_{\cou} & \to & F
& \to & 0 \\
  &  & \hookdownarrow &  & \hookdownarrow &
 & \downarrow & &  \\
 0 & \to & M_{V, \lin} & \to & V \otimes \oo_{\cou} &
\to & \lin & \to & 0
\end{array}
$$
where $W \hookrightarrow V$ is defined by 
$W^* := \mathrm{Im} (V^* \to H^0(N^*))$,
hence $W^*$ generates $N^*$, 
and we call
$F^*:= \ker (W^* \otimes \oo \twoheadrightarrow N^*)$.

Then $F$ is a vector bundle with no trivial summands.
Moreover the morphism $F \to \lin$ is not zero, as 
$W \otimes \oo$ does not map to $ M_{V, \lin}$.
We have to show that  $\rk F =1$ and that 
$W = H^0(F)$.
We distinguish the two cases $h^1(F) = 0$ or 
$h^1(F) \neq 0$.

$\bullet$ Let us suppose that $h^1(F) =0$.
Then 
$h^0(F) = \chi (F) = \rk F (\mu (F)+1-g)$.
On the other hand, 
 $h^0(F) > \rk F$ as
 $F$ is globally generated and not trivial.

Together this yields
\begin{equation}
\label{mug}
\mu(F)>g ~.
\end{equation}

Furthermore
\begin{equation}
\label{mintot}
\mu(N) = - \deg F / (\dim W - \rk F)  
\leqslant -\mu(F) / (\mu(F)-g) = \mu(M_F)~,
\end{equation}
as $\dim W \leqslant h^0(F) = \rk F (\mu (F)+1-g)$.

Consider the image $I = \mathrm{Im} (F \to \lin) \subseteq \lin $.
The commutative diagram 
$$\begin{array}{ccccc}
W & \hookrightarrow & H^0(F)& \to & H^0(I)\\
 \hookdownarrow &&&& \hookdownarrow \\
V &&\hookrightarrow && H^0(\lin)
\end{array}$$
shows that 
the map $W \to H^0(I)$ is injective and its image $W^{\prime} \subset H^0(I)$
is contained in 
$ V \subset H^0(\lin)$,
hence 
$N \hookrightarrow M_{W^{\prime},I} \hookrightarrow M_{V, \lin}$. 
As $N$ is a subbundle of $M_{V, \lin}$ of maximal slope, this yields
$\mu(N) \geqslant \mu(M_{W^{\prime},I})$,
\emph{i.e.}
 $-\deg F / \rk N \geqslant - \deg I / \rk M_{W^{\prime},I}$.
 Then
$$\deg F \leqslant \deg I (\rk N / \rk M_{W^{\prime},I} ) 
\leqslant \deg I \leqslant \deg \lin =d~.$$

If $\rk F \geqslant 2$, then 
$\mu(F) \leqslant \deg \lin / 2 = d/2$, 
so
$$\mu(N) \leqslant \frac{-\mu(F)}{\mu(F)-g} 
\leqslant  \frac{-d/2}{d/2 -g} = \frac{-d}{d-2g} \leqslant
\frac{-d}{d-g-c} = \mu (M_{V, \lin})  ~.$$
Here the first inequality is (\ref{mintot}). For the second
one shows that the function $ -x / (x-g)$ is strictly
increasing for $x>g$. Then use $\mu(F)>g$ due to (\ref{mug}).
Equality holds only if $\rk F =2$, $\deg F = d$, $W =H^0(F)$, and $g=c$. But in
this case we would find that $\dim W = h^0(F) = d+2 - 2g > d+1-g-c =
\dim V$, wich is impossible as by construction $W \hookrightarrow V$.

Hence  $\rk F =1$.
So $F = I$ is a globally generated and acyclic line bundle 
of degree $f \leqslant d$,
and $\mu(N) = -f/(\dim W -1)$.

It is easy to see that the case $f=d$ cannot hold, 
as in that case we cannot have 
$\mu(N) \geqslant \mu (M_{V, \lin})$. So $f\leqslant d-1$.

% If $ W \subsetneq H^0(F)$ then $\dim W -1 \leqslant f-g-1$.
% Hence
% $$ \mu(N) = -f/(\dim W -1) \leqslant -f / (f-g-1) 
% \leqslant -d/(d-g-1) = \mu (M_{V, \lin})~, $$ with equality only if
% $N = M_{V,\lin}$ (use that the function $ -x / (x-g-1)$ is strictly
% increasing for $x>g+1$, and that $\deg F>g+1 $).
% We indeed have $W = H^0(F)$.

% It remains to show that $f = \deg F \leqslant d -2$.
% Clearly, $f=d$ is impossible, 
% as it would lead to the absurd 
% $H^0(\lin) = W \subseteq V\subsetneq H^0(\lin)$.  
% Suppose $f = d-1$, then $h^0(F) = d-g = \dim V$,
% and  $W= H^0(F)\longiso V $.
% This is also impossible, as it would yield that 
% $V \cong H^0(F)$ gnerates the proper subsheaf $F$ but not $\lin$.
% So we have $f = \deg F \leqslant d-2$.

 $\bullet$ In the case $h^1(F) \neq 0$, by lemma \ref{butlem},
 $\mu(N)\leqslant -2$. Equality holds only if 
 $F = \omega_{\cou}$ and $W = H^0(\omega)$, 
 or if the curve $\cou$ is hyperelliptic and $F$ is 
 the hyperelliptic bundle. 
 In the latter case the only generating space of global
 sections is $H^0(F)$. 
 In any case we have $f = \deg F < d-1$.
 \end{proof}

\begin{rmk}

The diagram in the statement of the lemma is 
a construction from Butler's proof 
of theorem \ref{but}.

\end{rmk}

\begin{rmk}
\label{num}

Loking carefully at the numerical invariants in the above proof,
we can deduce some inequalities which will be useful in the following:
let us consider again the diagram in the above lemma
$$
\begin{array}{ccccccccc}
~~~ 0 & \to & N & \to & W \otimes \oo_{\cou} & \to & F
& \to & 0 \\
  &  & \hookdownarrow &  & \hookdownarrow &
 & \hookdownarrow & &  \\
 0 & \to & M_{V, \lin} & \to & V \otimes \oo_{\cou} &
\to & \lin & \to & 0 
\end{array}
$$
and suppose that 
 $h^1(F) = 0$. Let us call $f:= \deg F$, $s:= d-f$, and $b:= \cod_{H^0(F)} W $.
Then we can show that 
\begin{equation}
\label{bs}
0 < c-b < s \leqslant \frac{d}{g+c} (c-b) ~.
\end{equation}

In fact, as $W \hookrightarrow V$, and $W \neq V$, then 
$$d-s +1 -g -b = h^0(F) -b =  \dim W   < \dim V = d +1 -g -c   ~,$$
hence $c-b<s$. And as 
$$ -\frac{d-s}{d-s-g-b} = \mu (N)  \geqslant \mu( M_{V,\lin})= -\frac{d}{d-g-c}  ~,$$
then $s(g+c) \leqslant d(c-b)$, hence $c-b >0$ and 
$s \leqslant \frac{d}{g+c} (c-b)$.

\end{rmk}

\subsection{Line bundles of degree $d=2g+2$}
\label{2g+2}

A first consequence of these lemmas is the following proposition 
asserting semistability for hyperplane tranforms of line 
bundles of degree $2g+2$.

 \begin{prop}
 \label{2g+2semi}

 Let $\lin$ be a line bundle of degree $d=2g+2$  on a
  curve $\cou$ of genus $g \geqslant 2$.
 Then $M_{V, \lin}$ is semistable for every 
 generating
 hyperplane $V \subset H^0 (\lin)$.
 It is strictly semistable if 
 $\cou$ is hyperelliptic.
 
 \end{prop}

\begin{proof}
Let us prove the semistability of $M_{V, \lin}$.

Consider a stable subbundle 
$N \hookrightarrow M_{V, \lin}$ of maximal slope, and suppose that 
it destabilizes $M_{V, \lin}$ in the strict sense,
\emph{i.e.}
 $\mu(N) >- 2 = \mu (M_{V, \lin})$. 
By lemma \ref{part1} and remark \ref{num} (we have $b=0$ in this case),
we know that $N$ fits into a diagram
$$
\begin{array}{ccccccccc}
~~~ 0 & \to & N & \to & H^0(F) \otimes \oo_{\cou} & \to & F
& \to & 0 \\
  &  & \hookdownarrow &  & \hookdownarrow &
 & \hookdownarrow & &  \\
 0 & \to & M_{V, \lin} & \to & V \otimes \oo_{\cou} &
\to & \lin & \to & 0
\end{array}
$$
with $F$  a line bundle of degree $\deg F \leqslant d-2 = 2g$.
Moreover, $h^1(F) = 0$ since otherwise 
$\mu(N) \leqslant - 2$ by lemma \ref{butlem}.
Hence $\rk N = \deg F -g$, and 
$$\mu(N)= - \deg F / (\deg F -g) \leqslant -2g/ (2g -g) = -2 ~,$$ 
(again, use that the function $-x / (x-g)$ is strictly increasing for $x>g$).  
So it is not possible to find a strictly destabilizing $N$.

If the curve is hyperelliptic, then  $M_{V, \lin}$ is strictly semistable:
we can show that there is a line bundle of degree $-2$ injecting in 
 $M_{V, \lin}$.
In fact we can consider the line bundle $A$ 
dual of the only $g^1_2$ of the curve,
\emph{i.e} the dual of the hyperelliptic bundle.

The hyperelliptic bundle $A^*$ has $h^0(A^*)=2$, 
and
from the exact sequence 
$ 0 \to M_{V,\lin} \otimes A^*
  \to V \otimes A^*
\to \lin \otimes A^* \to 0 ~,$
we see that there are destabilizations of 
$M_{V, \lin}$ by the line bundle $A$ if and only if
$$H^0(M_{V,\lin} \otimes A^*)
= \ker (\varphi \colon V \otimes H^0(A^*)
\to H^0( \lin \otimes A^*)) \neq 0 ~. $$
%\begin{sloppypar}
Counting dimensions we see that 
the map $\varphi$ cannot be injective:

$$\dim V \cdot \dim  H^0(A^*)= (g+2) 2   > 
g+5 = \dim H^0( \lin \otimes A^*).$$
\end{proof}

In order to prove stability for non hyperelliptic curves though, we need to take a
generic
hyperplane, and not just a generating one.

The following is a special case of a more general result proven in section \ref{2g+2c}

 \begin{thm}
 \label{2g+2stab}

 Let $\lin$ be a line bundle 
 of degree $d=2g+2$  on a
  curve $\cou$ of genus $g \geqslant 2$.
 Then $M_{V, \lin}$ is stable for a 
 generic
 hyperplane $V \subset H^0 (\lin)$
  if and only if
 $\cou$ is non hyperelliptic.

 \end{thm}

\subsection{Line bundles of degree $d> 2g+2c$}

%\begin{rmk}

Here we show that for a generic 
subspace the transform of a line bundle of degree 
$d > 2g + 2c$ is stable.
In contrast to proposition \ref{2g+2semi},
we have to consider generic hyperplanes,
and not just generating ones.

%\end{rmk}

\begin{thm}
\label{princ}

Let $\lin$ be a line bundle  of degree $d$ on a curve 
$\cou$ of genus $g \geqslant 2$, such that  
$  d > 2g + 2c$, with $1 \leqslant c \leqslant g$.
Then $M_{V, \lin}$ is stable for a 
generic
subspace $V \subset H^0 (\lin)$ of codimension $c$.

\end{thm}

\begin{proof}

Let us proceed as in proposition \ref{2g+2semi}.
We have that $-2 < \mu (M_{V, \lin})< -1$.

Consider a stable subbundle 
$N \hookrightarrow M_{V, \lin}$ of maximal slope.
By lemma \ref{part1} we know it fits into a diagram
$$
\begin{array}{ccccccccc}
~~~ 0 & \to & N & \to & W \otimes \oo_{\cou} & \to & F
& \to & 0 \\
  &  & \hookdownarrow &  & \hookdownarrow &
 & \hookdownarrow & &  \\
 0 & \to & M_{V, \lin} & \to & V \otimes \oo_{\cou} &
\to & \lin & \to & 0 ~.
\end{array}
$$

We can right away conclude
that $h^1(F) = 0$,
as by  lemma
\ref{butlem} we would otherwise have 
$\mu(N) \leqslant -2$.

So $F$ is a globally generated line bundle with $h^1(F)=0$,
 $\deg F =: d-s \leqslant d-2$, and $W$ is a $b$-codimensional subspace 
of $H^0(F)$.
By remark \ref{num}, we see that for every $b$ with 
$0 \leqslant b <c$ there is a finite number of $s$ giving rise to 
a possible destabilization of $ M_{V, \lin}$.

For any of those $b$ and $s$  we will construct a parameter space
allowing $F$, $W$, and the subspace $V \subset H^0(\lin)$  to vary.

For any such $b$ and $s$ we want to consider the parameter space
$\mathcal{D}_{b,s}$,
parametrizing subspaces $V \subset H^0(\lin)$ together with 
a destabilizing bundle of  
$M_{V, \lin}$ of degree $s-d$ originating from a subspace $W$ as in
the construction above:
$$\mathcal{D}_{b,s} :=\{ (F~,~ F \hookrightarrow \lin ~,~ 
W \subset  H^0(F)) ~,~
V \subset  H^0(\lin)) ~|~ F \in \pic^{d-s} (\cou) ~,~ $$
$$(\varphi \colon F \hookrightarrow \lin) \in 
\pp(H^0(F^* \otimes \lin)) ~,~
W \in \gr  (b, H^0(F)) $$
$$ V \in \gr (c, H^0(\lin))~,~
 \varphi_{|W} \colon W \hookrightarrow V \subset H^0(\lin)  \} ~. $$

In order to estimate its dimension,
we use the natural morphisms 
$$\pi_{b,s} \colon \mathcal{D}_{b,s} \to \pic^{d-s} (\cou)~,~
(F, F \hookrightarrow \lin ,W, V ) \mapsto F~,$$
and
$\rho_{b,s} \colon \mathcal{D}_{b,s} \to \gr(c,H^0(\lin))$, 
$(F, F \hookrightarrow \lin ,W, V ) \mapsto V$.

The image of $\pi_{b,s}$ is formed by all the line bundles 
$F \in \pic^{d-s}(\cou)$ such that 
$h^0(F^* \otimes \lin) \ne 0$.
In particular 
$\dim \pi_s ( \mathcal{D}_s) = \min (s,g)$, 
because the degree of $F^* \otimes \lin$ is $s$.
The fiber over $F \in \pi_{b,s}(\mathcal{D}_s)$ 
has the same dimension as 
$\pp(H^0(F^* \otimes \lin)) \times \gr(b,(H^0(F))) \times
\gr(c,(H^0(\lin)/W))$.

By Clifford's theorem, $h^0(F^* \otimes \lin) = s/2 +1$ if
$s\leqslant 2g$,
and $h^0(F^* \otimes \lin) = s+1-g$ otherwise.
So, 
$$\dim \mathcal{D}_{b,s} \leqslant \min (s,g) 
+ \sup(s/2, s-g) + 
b(d-s-g+1-b) + c(s+b-c) \leqslant$$
$$ (3/2) s+ 
b(d-s-g+1-b) + c(s+b-c) ~.$$

{\bf Claim:} for $g,d,c$ as in the hypothesis and 
$s,b$ satisfying the inequalities of remark \ref{num}, 
we have
$$ (3/2) s+ b(d-s-g+1-b) + c(s+b-c) < c(d+1-g) - c^2 = 
\dim \gr (c, H^0(\lin))~.$$

% \begin{sloppypar}
Proving the claim, we show that for all $s$ and $b$ 
giving rise to possible destabilizations, the morphisms
\hbox{$\rho_{b,s} \colon \mathcal{D}_{b,s} \to  \gr(c,H^0(\lin))$} 
have a locally closed image of dimension strictly smaller than 
$\gr ( c, H^0(\lin))$, 
hence the generic subspace avoids all possible destabilizations of $M_{V, \lin}$.
% \end{sloppypar}

The claim is equivalent to
$$\frac{3s}{2(c-b)} +s +b < d+1-g ~,$$
using inequalities (\ref{bs}) we get
$$\frac{3s}{2(c-b)} +s +b \leqslant \frac{3/2 + (c-b)}{g+c}d +b ~,$$
hence we want to prove
$$ \frac{3/2 + (c-b)}{g+c}d +b < d+1-g ~,$$
which is equivalent to 
$$ \frac{b+g-1}{b+g- 3/2} < \frac{d}{g+c} ~,$$
and as $b\geqslant 0 \geqslant 2-g$ then 
$ \frac{b+g-1}{b+g- 3/2} \leqslant 2 < \frac{d}{g+c}$.
\end{proof}

\subsection{Line bundles of degree $d= 2g+2c$}
\label{2g+2c}

We have shown in section
\ref{2g+2}
that hyperplane tranforms of a degree ${2g+2}$ line bundle
are always semistable.

We prove here, that generic $c$-codimensional transforms 
of a degree ${2g+2c}$ line bundle
are stable, except in the hyperelliptic case, where they are strictly semistable.

\begin{thm}
\label{2g+2stab}

Let $\lin$ be a line bundle 
of degree $d=2g+2c$  on a
 curve $\cou$ of genus $g \geqslant 2$.
Then $M_{V, \lin}$ is semistable for a 
generic
subspace $V \subset H^0 (\lin)$ of codimension $c$.
It is stable
 if and only if
$\cou$ is non hyperelliptic.

\end{thm}

\begin{proof}
As in the proof of theorem \ref{princ}
we want to construct parameter spaces for destabilizations, and verify
by dimension count that the generic subspace avoids them.

Let us consider a line bundle $\lin$
of degree $d=2g+2c$  on a
 curve $\cou$ of genus $g \geqslant 2$, and the transform
$M_{V, \lin}$ for a 
subspace $V \subset H^0 (\lin)$ of codimension $c$.

To show semistability, let us suppose that there is a destabilizing
stable vector bundle $N \hookrightarrow M_{V, \lin}$, with 
$ \mu (N) > \mu (M_{V, \lin}) = -2$.

By lemma \ref{part1} we know it fits in the diagram
$$
\begin{array}{ccccccccc}
~~~ 0 & \to & N & \to & W \otimes \oo_{\cou} & \to & F
& \to & 0 \\
  &  & \hookdownarrow &  & \hookdownarrow &
 & \hookdownarrow & &  \\
 0 & \to & M_{V, \lin} & \to & V \otimes \oo_{\cou} &
\to & \lin & \to & 0 ~.
\end{array}
$$
and we can suppose that $h^1(F)=0$ by lemma 
\ref{butlem}.

In this case we can follow the same computations as in theorem \ref{princ}:
we have a parameter space for destabilizations 
$$\mathcal{D}_{b,s} :=\{ (F~,~ F \hookrightarrow \lin ~,~ 
W \subset  H^0(F)) ~,~
V \subset  H^0(\lin)) ~|~ F \in \pic^{d-s} (\cou) ~,~ $$
$$(\varphi \colon F \hookrightarrow \lin) \in 
\pp(H^0(F^* \otimes \lin)) ~,~
W \in \gr  (b, H^0(F)) $$
$$ V \in \gr (c, H^0(\lin))~,~
 \varphi_{|W} \colon W \hookrightarrow V \subset H^0(\lin)  \} ~,$$
whose dimension is bounded by
$$\dim \mathcal{D}_{b,s} \leqslant (3/2) s+ 
b(d-s-g+1-b) + c(s+b-c) ~,$$
with $b$ and $s$ satisfying 
$0 < c-b < s \leqslant \frac{d}{g+c} (c-b) $.

Except in the case $b=0$ and $g=2$,
 we can follow the very same proof of theorem \ref{princ}, and we see that this bound shows that the generic subspace avoids the destabilization locus.
 
 In the case $b=0$ and $g=2$ as well,  
 it can be easily shown that 
 $\dim  \mathcal{D}_{b,s}  < \dim \gr (c, H^0(\lin)) $, 
for all $s$ giving rise to destabilizations.

To show that we have strict semistability in the hyperelliptic case, we can proceed as in 
proposition \ref{2g+2semi}, and show that dual of  the hyperlliptic bundle
is a subbundle of $M_{V, \lin}$, of slope $-2$.

To show that we have stability in the non hyperelliptic case,
we have to exclude slope $-2$ subbundles $N \hookrightarrow M_{V, \lin}$.

Again we can apply lemma \ref{part1}
and consider the diagram
$$
\begin{array}{ccccccccc}
~~~ 0 & \to & N & \to & W \otimes \oo_{\cou} & \to & F
& \to & 0 \\
  &  & \hookdownarrow &  & \hookdownarrow &
 & \hookdownarrow & &  \\
 0 & \to & M_{V, \lin} & \to & V \otimes \oo_{\cou} &
\to & \lin & \to & 0 ~,
\end{array}
$$
where we can distinguish the two cases $H^1 (F) =0$, and $H^1 (F) \neq 0$.

In the case $H^1 (F) =0$ we can follow again the same computations as in theorem \ref{princ}.

In the case $H^1 (F) \neq 0$, lemma \ref{butlem} implies $F = \omega$ and $N= M_{\omega}$,
hence the parameter space for destabilizations will be
$$\mathcal{D} :=\{ (\omega \hookrightarrow \lin ~,~ 
V \subset  H^0(\lin)) ~|~  H^0 (\omega) \subset V \} ~,$$
and it can be shown that $\dim \mathcal{D} < \dim \gr (c, H^0 (\lin))$.

\end{proof}

\section{Theta divisors and transforms}

When a vector bundle has integer slope
$\mu(E) = \mu \in \zz$, we can define the set 
$$\Theta_E := \{ P \in \picard^{\nu}(C) ~|~ H^0(C, E\otimes P) \ne 0 \} ~,$$
where $\nu := g-1-\mu$.

As $\chi (E \otimes P) =0$, either $\Theta_E = \picard^{\nu}(C)$, or
it has a natural structure of
effective divisor in $\picard^{\nu}(C)$. 
In the latter case we say that
\emph{E admits a theta divisor}.
The class of this divisor in $H^2(\picard^{\nu}(C), \zz)$ is 
$\rk E \cdot \vartheta$, where $\vartheta$ is the class of the canonical
theta divisor of $\picard^{\nu}(C)$.

Whenever a vector bundle admits a theta divisor, then it is semistable.
And strictly semistable vector bundles admitting a theta divisor have 
\emph{non integral} theta divisors.

However, there are examples of stable vector bundles with no theta divisor,
or with a reducible theta divisor.

Beauville shows in \cite{bove} that the total transform 
$M_L$ of a degree $2g$ line bundle $L$ on a genus  $g$ curve $C$
always has a reducible theta divisor.
And that if 
$L$ is very ample, and $C$ is not hyperelliptic, then $M_L$ is stable.

The vector bundles considered above, 
\emph{i.e.} tranforms of degree $d \geqslant 2g + 2c$ line bundles,
with respect to $c$ codimensional subspaces of global sections,
have slope $\mu$ such that $-2 \leqslant \mu < -1$.
The case of integer slope $\mu = -2$ appears if and only if $d = 2g + 2c$.

Following the same argument as in \cite{bove}, we prove that
for the generic $V \subset H^0(C, \lin)$ 
those tranforms always carry a non integral theta divisor.

To prove that, for a generic $V \subset H^0(C, \lin)$ within 
the numerical conditions above,
the transform $M_{V, \lin}$ admits a theta divisor,
we need the following lemma:

\begin{lem}
\label{geninter}

Let $P$ be a $2$-dimensional vector space, 
$H$ a vector space of dimension $n+c$, and 
$K \subset P \otimes H$ a subspace of dimension $2c$.
If $K$ contains no pure vectors, then 
the generic $n$-dimensional subspace $V \in \gr (n, H)$ verifies
$$ K \cap (P \otimes V) = 0~.$$

\end{lem}

\begin{proof}
We consider the map
$$
\begin{array}{ccc}
f \colon \gr(n, H) & \to & \gr (2n, P \otimes H)\\
V & \mapsto & P \otimes V ~,
\end{array}
$$
and we claim that the image of $f$ is not contained 
in the closed subscheme 
$$Z := \{ W \in \gr (2n, P \otimes H) ~|~ \dim K \cap W \geqslant 1\}~.$$
Let us observe at first that $Z$ carries a filtration
$$Z = Z_1 \supseteq Z_2 \supseteq \dots \supseteq
Z_s := \{ W \subset P \otimes H ~|~ \dim K \cap W \geqslant s\}
\supseteq \dots~.$$
The tangent space of the grassmannian $\gr (2n, P \otimes H)$ at a point
$W$ is 
$$T_W \gr (2n, P \otimes H) = \Hom (W , P \otimes H / W) ~.$$
The subscheme $Z_s \setminus Z_{s+1}$ is smooth and its tangent space
at a point $W$, is given by  first order deformations of 
$W \subset P \otimes H$ 
that deform $W\cap K$ into an $s$-dimensional subspace of $K$:
$$T_W (Z_s \setminus Z_{s+1}) = 
\{ \varphi \in \Hom (W, P \otimes H/ W)
~|~ \varphi(W \cap K) \subseteq K/ (W\cap K) \}~.$$
And the differential of the morphism $f$ at the point $V\in {\gr(n, H)}$
 is the map
$$
\begin{array}{cccc}
\de f_V \colon & T_V {\gr(n, H)} & \to & T_{P \otimes V} \gr (2n, P \otimes H)\\
& \varphi \in \Hom (V, H/V)
& \mapsto & 1 \otimes \varphi \in \Hom (P \otimes V, P \otimes (H/V)) ~.
\end{array}
$$

We can prove now that if $V \in {\gr(n, H)}$ is a subspace
such that  $P \otimes V \in Z_s \setminus Z_{s+1}$,
then $\de f_V(T_V {\gr(n, H)} ) \nsubseteq T_{P \otimes V} (Z_s \setminus Z_{s+1})$:
we claim that there exists a $\varphi \in  \Hom (V, H/V)$
such that $1 \otimes \varphi \not \in T_{P \otimes V} (Z_s \setminus Z_{s+1}) $.

To see this, let us choose a basis $(e_1 , e_2)$ for $P$,
and a vector
$w = e_1 \otimes v_1 + e_2 \otimes v_2 \in K \cap (P \otimes V)$.
By the hypothesis on $K$, $v_1 \nparallel v_2$.
Let us consider now a vector 
$z =  e_1 \otimes z_1 + e_2 \otimes z_2 \in 
(P \otimes (H/V))$ such that $z \notin (K/(P\otimes V \cap K))$.

Then if we choose a 
$\varphi \in \Hom (V, H/V)$ such that 
$\varphi(v_1) = z_1$ and $\varphi (v_2) = z_2$,
we have that 
$(1 \otimes \varphi) (w) = z \notin K/(P\otimes V \cap K)$.
Hence the image of a generic deformation of $V$ avoids the subscheme 
$Z \subset \gr (2n, P \otimes H)$.

\end{proof}

We can now prove the existence of theta divisors for generic transforms 
of slope $-2$.

\begin{thm}

Let $\lin$ be a line bundle of degree $d = 2g +2c$ on a 
genus $g$ curve $C$, where $c \in \nn$ is a positive integer and $g \geqslant 2$.
Then, if $V \subset H^0(C, \lin)$ is a generic $c$-codimensional 
subspace, the transform $M_{V, \lin}$ admits a non integral theta divisor.

\end{thm}

\begin{proof}
We recall that $\mu(M_{V, \lin}) = -2$.
We have to show first that, for the generic $V \subset H^0(C, \lin)$,
$\Theta_{M_{V, \lin}} \neq \picard^{g+1} (C)$, \emph{i.e.} that
there is a $P \in \picard^{g+1} (C)$ such that 
$H^0(M_{V, \lin} \otimes P) = 0$.
By the exact sequence  
$$ 0 \to M_{V, \lin} \otimes P \to V \otimes P \to \lin \otimes P \to 0$$
this is the same as a $P \in \picard^{g+1} (C)$ such that
the multiplication map
$$\mu \colon V \otimes H^0(P) \to H^0(\lin \otimes P)$$ is injective.

If $P$ belongs to the divisor $D = (\omega_C) - C_{g-2} +C \subset \picard^{g+1} (C)$, 
\emph{i.e.} if $P$ can be written in the form
$P = \omega_C (x_1 -x_2 - \dots - x_{g-1})$ for some points
$x_1, x_2, \dots, x_{g-1} \in C$,
then either $h^0(P) >2$, or $h^0(P) =2$ and $P$ has a base point.
In any case this implies that $\mu$ is not injective for any $V$
(c.f. \cite{bove}).

Any $P$ in $\picard^{g+1} (C) \setminus D$ is base point free and has $h^0(P)=2$.
Let us fix such a $P$, and assume by generality that $h^1(L \otimes P^*)=0$.
We claim that for the generic $V \subset H^0(C, \lin)$ of codimension $c$ 
the multiplication map $\mu \colon  H^0(P) \otimes V \to H^0(P \otimes \lin)$
is injective.
From the exact sequence
$$ 0 \to P^* \to H^0( P) \otimes \oo_C \to \ P \to 0$$
we get
$$ 0 \to H^0(P^* \otimes { \lin} ) \to H^0(P)  \otimes H^0(\lin) 
\to H^0(P \otimes \lin) \to 0~,$$
hence the map $\mu$ is injective if and only if 
the subspace $V \subset H^0(C, \lin)$ verifies
$ H^0(P^* \otimes { \lin} ) \cap (H^0(P) \otimes V) =0$.
And this is given by lemma \ref{geninter}.

Hence we know that for the generic subspace $V$, the transform 
$M_{V, \lin}$ admits a theta divisor.
To observe that it is not integral,
we notice that the set of points of $\Theta_{M_{V, \lin}}$ contains
the divisor $D$ whose cohomology class is $(g-1) \vartheta$
(cf. \cite{fmp}).
As the cohomology class of $\Theta_{M_{V, \lin}}$ is 
$(g+c) \vartheta$, it must be a non integral divisor.

\end{proof}

As we have proved the existence of theta divisors for transforms with respect 
to subspaces of any codimension, 
this shows semistability in some cases not previously treated:

\begin{cor}

Let $\lin$ be a line bundle of degree $d = 2g +2c$ on a 
genus $g$ curve $C$, where $c \in \nn$ is \emph{any} positive integer and $g \geqslant 2$.
Then, if $V \subset H^0(C, \lin)$ is a generic $c$-codimensional 
subspace, the transform $M_{V, \lin}$ is semistable.

\end{cor}

\begin{rmk}

If $C$ is not hyperelliptic and $\lin$ is a degree $d = 2g + 2c$ line bundle,
where $d \notin 2(g-1)\nn$,
then 
the transform $M_{V, \lin}$ of $\lin$ with respect to 
a generic subspace $V \subset H^0(C, \lin)$ of codimension $c$
admits a reducible  theta divisor.
In fact the set of points 
 of $\Theta_{M_{V, \lin}}$ contains
the divisor 
$$D= (\omega_C) - C_{g-2} +C \subset \picard^{g+1} (C)~,$$
which is irreducible if $C$ is not hyperelliptic,
and whose cohomology class is $(g-1) \vartheta$.
As  the cohomology class of $\Theta_{M_{V, \lin}}$ 
is $(g+c) \vartheta = (d/2) \vartheta$, it cannot be a multiple
of $(g-1) \vartheta$, then $\Theta_{M_{V, \lin}}$ must be reduciblle.

Hence, if $c \leqslant g$ and $c \neq g -2$ we have
further examples of stable 
vector bundles (by theorem \ref{2g+2stab})
with reducible theta divisors.

\end{rmk}

%\begin{rmk}

%By reducible theta divisor we mean that the natural scheme structure 
%of $\Theta_E \subset \picard^{\nu} (C)$ is not integral.
%An example of a theta divisor which should be considered reducible,  
%even though its scheme structure is irreducible,
%is the theta divisor of the strictly semistable vector bundle 
%$E = L \oplus L$, where $L$ is a line bundle on $C$, 
%and $\Theta_E = 2 \Theta_L$.

%\end{rmk}

\section{Conclusions}

We have proven stability of transforms of line bundles 
with respect to subspaces 
of low codimension. On the converse, it is rather easy to show the stability 
of transforms with respect to subspaces of low dimension:
any stable vector bundle $M^*$ of slope $\mu (M^*) > 2g -1$ 
is  globally generated.
Hence we can pick any stable vector bundle $M^*$ 
of determinant $\lin$ and rank $r$, such that
$r < {d}/{(2g-1)}$, where $\deg \lin = d$.
Choosing any generating subspace $V^* \subset H^0(M^*)$ of rank $r+1$, 
we get an exact sequence
$$0 \to \lin^* \to V^* \otimes \oo \to M^* \to 0  ~.$$
Dualizing we get an exact sequence 
$$0 \to M \to V \otimes \oo \to \lin \to 0  ~,$$
where M is a stable transform of $\lin$.
Hence, every stable bundle of rank $r < {d}/ {(2g-1)}$ and determinant
$\lin^*$, is a stable transform of $\lin$. So the rational map 
$\gr (r+1, H^0(\lin)) \dashrightarrow \mathcal{SU}(r,\lin) $ is dominant.

By the same argument 
we see that there is only one globally generated vector bundle, 
among vector bundles  of determinant $\lin$ and rank $d-g$ with no trivial summands,
where $d= \deg \lin \geqslant 2g$.
Furthermore this is semistable, and
even stable if $d>2g$.
In fact having such a globally generated bundle $N$,
we can pick a vector space $V$ of global sections of dimension $\rk N +1$ generating $N$.
This gives rise to the exact sequence
$$0 \to \lin^* \to V \otimes \oo \to N \to 0  ~,$$
and dualizing
$$0 \to N^* \to V^* \otimes \oo \to \lin \to 0  ~.$$

But as $N$ is globally generated
and has no trivial summands, then $H^0(N^*)=0$.
And since $V^*$ and $H^0(\lin)$ have the same dimension, 
then $V^* \iso H^0(\lin)$.
Hence $N^* = M_{\lin}$ is unique.

So when we consider the rational map 
$\gr (r+1, H^0(\lin)) \dashrightarrow \mathcal{SU}(r,\lin) ~,~  
V \mapsto (M_{V, \lin})^*$,
we are saying that its image is made by globally generated bundles,
and we can sum all this up in the following table, where we suppose that 
$d>2g +2c$, with $1 \leqslant c \leqslant g$:

$$
\begin{tabular}{|c|c|c|}
\hline
$\mathrm{rk} (M_{V, \lin}) =r$ &   stability & map \\
\hline
$1 \leqslant r < {d}/({2g-1})$  &
 stable & $\gr (r+1, H^0(\lin)) \dashrightarrow \mathcal{SU}(r,\lin) $
 dominant\\
 \hline
  $\frac{d}{2g-1} \leqslant r < d-g-c $ &
  \textbf{??} & \textbf{??} \\
  \hline
  $d-g-c \leqslant r < d-g$ &
  stable
&$\gr (r+1, H^0(\lin)) \dashrightarrow \mathcal{SU}(r,\lin) $\\
\hline
$r=d-g$   & stable &
$ \{*\} \hookrightarrow \mathcal{SU}(r,\lin) $\\
\hline
\end{tabular}
$$
where theorem \ref{princ}
corresponds to the existence of the rational map
$$\gr (r+1, H^0(\lin)) \dashrightarrow \mathcal{SU}(r,\lin) ~.$$

\bibliographystyle{amstronzo}
\bibliography{bibliostabtrans}

 \vspace{1cm}

 \begin{flushright}

 \small{

 Ernesto Carlo \textsc{Mistretta}

 \texttt{ernesto@math.jussieu.fr}

 \textsc{Institut de Mathématiques de Jussieu}

 Équipe de Topologie et Géométrie Algébriques

 175, rue de Chevaleret

 75013 Paris}

 \end{flushright}

\end{document}